\documentclass[14pt]{extarticle}
\usepackage[margin=1in]{geometry}
\usepackage{dsfont}
\usepackage{amsfonts}
\usepackage{tikz}
\usetikzlibrary{matrix}

\newtheorem{theorem}{Theorem}[section]
\newtheorem{pro}{Proposition}[section]
\newtheorem{lemma}{Lemma}[section]

\newtheorem{remark}{Remark}[section]
\newtheorem{cor}{Corollary}[section]

\newcommand{\proof}[1]{\noindent{\it\bf Proof:#1\ }}
\newcommand{\QED}{\hfill$\Box$\medskip}

\begin{document}

\title{ Higher-degree Smoothness of Perturbations I}
\author{Gang Liu }
\date{May 31,  2019}
\maketitle

\section{Introduction}
In \cite{5},  we developed the global perturbation method for constructing virtual moduli cycles for general symplectic manifolds. It was used to extend the definitions of Gromov-Witten invariants and Floer homology 
from the case of semi-positive symplectic manifolds to general case. This method is in the setting of infinite dimensions and is relatively simple (comparing with the other infinite dimensional methods in \cite{1}).

However, there was an analytic difficulty in this method. Recall that the global perturbations are constructed as follows. First the local perturbations on local uniformizers/slices $S_f$ were constructed and by multiplying smooth cut-off functions on the local slices, these local perturbations are globalized. Here $S_f$ is the local slice of the local chart ${\tilde W}(f)$ under the (local) group action of $G$, where   $G$ is the group of the reparametrizations of the domains, and  
  ${\tilde W}(f)$  consists of stable $L_k^p$-maps near the smooth center $f$.  
   Since the action of $G$ on  ${\tilde W}(f)$ is only continuous not even of class $C^1$, 
   the  transition functions between different local
  slices are only  of class $C^0$. Consequently, the 
  smooth local perturbations and  cut-off functions on one slice may not be smooth anymore viewed in the other slices (see a discussion for this in   \cite {6} ).  
  
  To overcome this difficulty, for the genus zero case with fixed domain $S^2$ and $G\simeq PSL(2, {\bf C})$, several methods were developed  
 \cite{8, 12, 13}. All these methods start with a smooth section $\xi=\xi_{S_f}$ with certain special properties on a local slice $S_f$,  end with   proving  that its $G$-equivariant extension $\xi_{{\cal O}_{S_f}}$ on the $G$-orbit ${{\cal O}_{S_f}}$ is smooth as well,  utilizing the special properties on the initial section $\xi$.

In this  paper, we  approach the problem of the lack of smoothness 
above with rather different  perspectives.

(i) We give necessary and sufficient conditions on  when the $G$-equivariant extension  is smooth, instead of only focusing on proving the smoothness of the 
$G$-equivariant extension for the special sections used in Gromov-Witten and Floer theory. 
(ii) We reproduce the simple conceptual proof of the smoothness of the cut-off function and its $G$-equivariant extension in \cite {13} 
 (compare the proofs in \cite {8, 9,  10, 12}), and give a new proof of the $C^{m_0}$-smoothness of the total evaluation map  as a corollary of Theorem 1.1 (compare the proof in \cite {11}).
 
  (iii) Above all we  find a  new point view on {\it lack of smoothness of G-action} that makes this negative  aspect into a  positive and crucial analytic input
   stated in Theorem 1.1, from which  all  analytic results  used in GW/Floer theory related to the  lack of smoothness, such as the results in (i) and (ii) above can be treated in a uniform and simple manner. 
   
   (iv) In such a treatment, the proofs should be given in such a way that  they  can be easily generalized.
   
   It is our hope that   our treatment here paves the way to 
   generalize part of the work of global analysis in \cite{7} for mapping spaces by  allowing  deformations of domains with various 
   topological types as well as   
   actions of the  groups of reparametrizations.   
  
  We now  explain the   new point view that  this work is based upon.
  
The prototype of the usual lack of smoothness  is the 
 statement that  the translation  action map $\Psi=\Psi_0:G\times L_k^p\rightarrow  L_k^p$ with $G={\bf R}^1$ and  $L_k^p=L_k^p({\bf R}^1)$ is only continuous
but not uniformly continuous and certainly  not $C^1$.  We refer to this statement as $S_0$. 

Statement $S_0$, a simple analytic fact  has a nature generalization  proved in \cite {14} 
 that    $\Psi_m:G\times L_k^p\rightarrow  L_{k-m}^p$ is of class $C^m$ if $m\leq k$. We refer to this generalization as statement $S_m $.

Statement  $S_m $, the more complete and precise version of the  {\it smoothness/lack of smoothness} of the action maps $\Psi_m$, in more general context, can be viewed in two opposite ways.

In the  negative direction, the usual version of the lack of the smoothness  stated in $S_0$ leads to the well-known analytic difficulties in the infinite dimensional method in the foundation of GW theory.
 
On the other hand, statement $S_m$ above is the prototype of the key analytic input in this paper on the $C^m$-smooth
of the  action maps $\Phi_m$ stated below, which naturally leads to the above mentioned criterion on the existence of smooth equivariant extensions (i). 
Combing this  with the simple proof of the smoothness
 of cut-off functions, the above analytic difficulties can be resolved  in a simple and direct  manner.

\begin{theorem}\label{Act1}
	Consider the action map  $\Phi:G_e\times L_k^p\rightarrow  L_k^p$  induced by the map $\phi$ below,  defined by $\Phi(\gamma, \xi)=\xi\circ \gamma, $ where $ L_k^p=:L_k^p(\Sigma, {\bf R}^N). $
	Its  restriction  	$\Phi_m:G_e\times L_{k+m}^p\rightarrow  { L}_{k}^p$ is of class $C^m$. Here we only assume that the "reparametrization" map $\phi:G_e\times \Sigma\rightarrow \Sigma$ is a  smooth  family of diffeomorphisms  of class $C^{\infty}$ parametrized by  $G_e\simeq B^n_{\epsilon}$ of a $\epsilon$-ball in ${\bf R}^n$, not necessarily  a (local) group action.

\end{theorem}

 We now explain how above theorem natural leads the necessary and sufficient conditions for the smoothness of $G_e$-equivariant extensions in the particular case that the ambient space  is  $L_k^p=:L_k^p(\Sigma, {\bf R}^N)$. %In this way, the complete form of the lack of smoothness resolves this main analytic difficulty by itself.
 
 To this end, let ${\tilde W}={\tilde W}_{k, p}$ be a small neighborhood of a smooth $f\in L_k^p$ above and $S_f=(S_f)_{k, p}$ be a local slice of the $G_e$-action. We may assume that ${\tilde W}$ is the same as the $G_e$-orbit ${\cal O}_{S_f}$.  Then
 for any $k\in {\tilde W}$, there  exists a
 $T^{-1}(k)\in G_e$ such that $k\circ T^{-1}(k)\in S_f.$ It was proved in sec. 3 of this paper and \cite{12} that $T:{\tilde W }\rightarrow G_e$ is of class $C^{m_0}$ with $m_0=[k-2/p]$. 
 
 Let $\eta: S_f\rightarrow L_k^p$ be a smooth map.
 Clearly  the  $G_e$-equivariant extension $\eta_{{\cal O}_{S_f}}$ is given by $\eta_{{\cal O}_{S_f}}(k)={\eta}_{S_f} (k\circ T^{-1}(k))\circ T(k)$ for $k\in W(f, {\bf H})$. It is decomposed as four maps 
 
 \begin{eqnarray*}
 	k&\rightarrow&(T^{-1}(k), k)\\
 	&\rightarrow&(T(k), k\circ T^{-1}(k))\\
 	&\rightarrow&(T(k), \eta(k\circ T^{-1}(k)))\\
 	&\rightarrow&(\eta(k\circ T^{-1}(k))\circ T(k)).
 \end{eqnarray*}

 Hence  
 $\eta_{{\cal O}_{S_f}}=\Phi^L_{m}\circ (Id_{G_e}\times   \eta_{-m, m})\circ  (IN_{G_e}\times \Phi^W_m)\circ (T^{-1}, Id^W)$.   The first and the second maps $(T^{-1}, Id^W):{\tilde W}_{k, p}\rightarrow G_e\times {\tilde W}_{k, p}(f)$ and $IN_{G_e}\times \Phi^W_m :G_e\times{\tilde W}_{k, p}(f)\rightarrow G_e\times({ S_f})_{k-m, p}$ are of class $C^{m}$ by the  Theorem 1.1, where $IN_{G_e}:G_e\rightarrow G_e$ is the inverse map. Note that  the map $\Phi^W_m :G_e\times{\tilde W}_{k, p}(f)\rightarrow {\tilde W}_{k-m, p}(f)$ should have the image ${\tilde W}_{k, p}(f)$ in general, but for the case here  for the elements $(T^{-1}(k), k)$, the corresponding image  is inside the submanifold $({ S_f})_{k-m, p}$.
 The last map $\Phi^L_{m}:
 G_e\times L_{k+m}^p\rightarrow L_{k}^p$ is of class $m$ by the Theorem 1.1 again. The third map  $Id_{G_e}\times \eta_{-m, m}:G_e\times({ S_f})_{k-m, p}(f)\rightarrow G_e\times L_{k+m}^p$   is of class $C^{m}$ only when we require that the map $\eta:(S_f)_{k, p}\rightarrow L_k^p$ can be extended into a $C^{m}$-map $\eta_{-m, m}:({ S_f})_{k-m, p}(f)\rightarrow L_{k+m}^p$ with image lying in $L_{k+m}^p\subset L_k^p.$
 
This gives the necessary and sufficient   condition  for the $C^m$-smoothness of  $\eta_{{\cal O}_{S_f}}$ in the above particular case. The proof    in next section is to deduce the general case to the special case  here.

To state the   necessary and sufficient   condition for the  case used in Gromov-Witten theory, we recall  some  notations in \cite {5} first.
Let $(M, \omega, J)$ be a compact smooth symplectic manifold with a $\omega$-compatible almost complex structure $J$. It gives rise a Riemannian metric $g_J(-, -)=\omega (-, J-)$ on $M$.  Denote the space of $L_k^p$-maps $f:\Sigma=S^2 \rightarrow M$ by  ${\widetilde {\cal B}}=:{\widetilde {\cal B}}_{k, p}$. 
 The  group $G= PSL(2, {\bf C})$ acts on  ${\tilde B}$ as the group of reparametrizations. Let  ${\widetilde {\cal B}}^s$ be the collection of  stable $L_k^p$-maps in the sense of \cite{16}.
Note that  the action of  $G$  on  ${\widetilde{\cal B} }^s$  is locally free in $G$. Since for a {\bf fix} element in $G$,  its  action is a smooth automorphism  of   ${\widetilde{\cal B} }^s$ and  related structures, for the purpose here, it is sufficient  to only consider local
action where $G$ is replaced by a small neighborhood $G_e$ of the identity (the notation used in the above theorem).

For $f\in {\widetilde {\cal B}}^s$,
let  $W(f, {\bf H})$ (defined in Sec. 3) be a local uniformizer inside a local chart ${\tilde W}(f)$ of ${\widetilde {\cal B}}^s$. Then $W(f, {\bf H})$ is a (local) slice of the action $G_e$ on ${\tilde W}(f)$, denoted also by $S_f$. We may assume that its $G_e$-orbit ${\cal O}_{S_f}={\tilde W}(f)$. 

Throughout  this paper we assume that (1) $m_0=[k-2/p]>0$  where $ k-2/p$ is the Sobolev smoothness of $L_k^p$-maps. Note that under this assumption, the space of $L_k^p$-functions on $\Sigma$ is a Banach algebra; (2) the center $f$ of each local chart is of class $C^{\infty}$ or sufficiently smooth.

Let ${\cal L}=:{\cal L }_{k,p}\rightarrow {\widetilde {\cal B}}$ be the Banach bundle with the fiber ${\cal L}_h=L_k^p(\Sigma, \wedge^1_{\Sigma}(h^*TM))$ for $h\in {\widetilde {\cal B}}$.

 For any $h\in  {\widetilde W}(f)$ and $\xi\in { {\cal L}}_f$,  the local trivialization ${\Pi}_f:{\widetilde W}(f)\times { {\cal L}}_f\rightarrow 
{ {\cal L}}|_{{\widetilde W}(f)}$ is given by ${\Pi}_f(h, \xi)(x)=P_{h(x)f(x)}\xi(x)$ for $x\in S^2$  where  $ P_{h(x)f(x)}$ is the map induced   by  the $J$-invariant parallel transport from $f(x)$ to $h(x)$ along the shortest
 geodesic.
  The restriction of  this local trivialization  to  ${ W}(f, {\bf  H})$ gives a $C^{\infty}$-smooth
local trivialization of the   bundle ${\cal L}\rightarrow { W}(f, {\bf  H})$.

 For a section $\eta:S_f\rightarrow {\cal L}_{k, p}|_{S_f}$ of class $C^{m_0}$,    
 with respect to the  local trivialization ${\cal L}|_{S_f}\simeq S_f\times L_k^p(\Sigma, \wedge^1_{\Sigma}(f^*TM))$,  $\eta(h)=(h, [\eta](h))$. Here $[\eta]:S_f\rightarrow  {\cal L}_{f}=L_k^p(\Sigma, \wedge^1_{\Sigma}(f^*TM))$ is the corresponding 
 $C^{m_0}$-function. 
 %Here 
 %$[\eta](h)=:P_{f, h}(\eta(h))$ and $P_{f, h}$ is the parallel transport from $h$ to $f$ along (the family of) shortest geodesics.
 
 Now we impose the conditions $C_1$ and $C_2$ on $[\eta]$ /$\eta$  defined as  follows. %Denote $S_f$ by $(S_f)_{0}$ consisting of $L_k^p$-maps. Then the corresponding space of $L_{k\pm m}^p$-maps will be denoted by $(S_f)_{\pm}$. Similarly 
 
 $C_1(=C_1(m)):$ The section   $\eta:S_f=(S_f)_{0}\rightarrow  {\cal L}|_{S_f}=({\cal L}|_{S_f})_0$ can be extended into  a $C^{m}$-smooth section
 $\eta_{-m}:(S_f)_{-m}\rightarrow  ({\cal L}|_{S_f})_0$ for some non-negative integer $m\leq m_0$. Here $(S_f)_{-m}=:W_{k-m, p}(f,{\bf H})$ consists of the $L_{k-m}^p$-maps.
 
 $C_2(=C_2(m)):$ For $h\in (S_f)_{-m}, $ the image of  $\eta_{-m}(h)$ is  lying in the fiber  of $L_{k+m}^p$-sections
 in the sense that $[\eta_{-m}](h)\in L_{k+m}^p(\Sigma, \wedge^1_{\Sigma}(f^*(TM)))$  for some  $m\leq m_0$.
 
  Note that here we have used the fact that the bundle $f^*(TM)\rightarrow \Sigma$ is of class $C^{\infty}$ or sufficiently smooth so that $L_{k+m}^p(\Sigma, \wedge^1_{\Sigma}(f^*(TM)))$ is well-defined.
 
  Hence   the  map $[\eta_{-m, m}]:(S_f)_{-m}\rightarrow L_{k+m}^p(\Sigma, \wedge^1_{\Sigma}(f^*(TM)))$  and the corresponding  section $\eta_{-m, m}$ are   of class $C^m$.
 
 \begin{theorem}\label{SmExt}
 	Under the  conditions $C_1$ and $C_2$  on  the smooth section $\eta=:\eta_{S_f}$ defined on the local slice $S_f$, the $G$-equivariant extension $\eta_{{\cal O}_{S_f}}$	 is at least of class $C^{m}$ for any $m\leq m_0$. 
 	Moreover, these conditions are necessary.
 \end{theorem}
 \medskip
 \medskip
 \medskip
 
% \begin{definition}
 %	A section $\xi$ of ${\cal L}\rightarrow S_f$ is  said to be  $\Pi$-admissible if is  satisfies   $C_1+C_2$ with respect to a given "standard" trivialization $\Pi=\Pi_{\nabla}:{\cal L}\simeq  S_f\times {\cal L}_f$  induced by a  connection $\nabla$ on $TM\rightarrow M$.

% \end{definition}

 \medskip
 \begin{remark}
 	Let  $\eta_{-m, n}:
 	(S_f)_{k-m, p}\rightarrow {\cal L}_{k+n, p}$ be a   smooth   section of class $C^l$.  It  gives rise  a section $\eta=\eta_{0, 0}:
 	(S_f)_{k, p}\rightarrow {\mathcal L}_{k, p}$ by first  restricting  $\eta_{-m, n}$ to $ (S_f)_{k, p}$ then composing  with the obvious inclusion map. Then $\eta$  satisfies the condition $C_1(m)$ and $C_2(n)$.    In the case that $\eta_{-m, n}$ is a bounded linear  operator   between  the corresponding 
 	Sobolev spaces,   this is  exactly the so called smoothing operator of degree $m+n.$ As for our case, since   $\eta_{-m, n}$ may not be   linear, the  smoothness of  $\eta_{-m, n}$ is imposed as a independent condition.  
 	
 	Thus   the  sections   satisfying the conditions $C_1(m) +C_2(n)$   can be identified  exactly with sufficiently smooth ( possibly non-linear) smoothing operators of degree $m+n.$ 
 	There are sufficient such sections in order to achieve transversality.   In fact for the applications in Gromov-Witten and Floer  theories, it is sufficient only considering  the constant sections in \cite{5}.
 \end{remark}
 
 Now recall the definition of "constant"  sections ${\xi_{S_f}}$  in  \cite {12}.
 
 Given a  $C^{\infty}$  element  $\xi_0 \in L_k^p(\Sigma, \wedge^1_{\Sigma} (f^*TM
 ))$ of the central fiber, let $[{\xi_{S_f}}]: S_f\rightarrow L_k^p(\Sigma, \wedge^1_{\Sigma} (f^*TM
 )) $ be the  constant function defined by $[{\xi_{S_f}}](h)=\xi_0$. Then  the "constant"  section   ${\xi_{S_f}}$  is  defined to be the corresponding  section under  the  local trivialization 
  ${\cal L}|_{S_f}\simeq S_f\times L_k^p(\Sigma, \wedge^1_{\Sigma}(f^*TM))$.  In   \cite {5, 8}, $\xi_0$ is obtained from  the elements in the   cokernel $K_f$ of the linearization at $f$ of the ${\overline {\partial}}_J$-operator in GW/Floer theory. These particular  constant sections are important for  GW/Floer theory as they are used to achieve the transversality for the perturbed ${\overline {\partial}}_J$-operator.

It is easy to see that the conditions $C_1$ and $C_2$  are satisfied for ${\xi_{S_f}}$. Indeed, the extended map   $[(\xi_{S_f})_{-m, m}]:(S_f)_{-m}\rightarrow L_{k+m}^p(\Sigma, \wedge^1_{\Sigma}(f^*(TM)))$  corresponding  to the section $(\xi_{S_f})_{-m, m}$ is still the constant map, $[(\xi_{S_f})_{-m, m}](h)=\xi_0\in C^{\infty}(\Sigma, \wedge^1_{\Sigma}(f^*(TM)))$ for any $h\in (S_f)_{-m}$.  

We refer the readers to \cite{ 13} for the definition of the geometric sections ${\eta}_{S_f}$. These sections are  obtained by pulling-back a "constant"-section of $TM\rightarrow M$ over a small ball $D(x_0)\subset M$ centered at $x_0$. Hence ${\eta}_{S_f}$ has similar properties as ${\xi_{S_f}}$. In particular, one can easily  show that  
the conditions $C_1$ and $C_2$  are  satisfied.

The following corollary recovers the results  in \cite {8} and  \cite{ 13}.

 \begin{cor}
 	Both the constant  sections ${\xi_{S_f}}$  and  the geometric sections ${\eta}_{S_f}$   are of class $C^{m_0}$  viewed in any   local slices. % Moreover ${\xi'_{S_{f'}}}$ above sill satisfies the two conditions $C_1$ and $C_2$ although it is not a "constant" section anymore. 
 \end{cor}
 %More generally for each  $\phi\in G_e$ and the corresponding local slice $S_{f\circ \phi}$, let ${\cal L}|_{S_f}={\cal L}|_{S_{f\circ \phi}}=\simeq S_{f\circ \phi}\times L_k^p(\Sigma, \wedge^1_{\Sigma}((f\circ \phi )^*TM))$ be the  local trivialization of class $C^{m_0}$. Here $L_k^p(\Sigma, \wedge^1_{\Sigma}, ((f\circ \phi )^*TM)))$ is the central fiber.
 Now assume that $p$ is an even integer.
 It was proved in last section and \cite {8, 12, 13} that   both the $p$-th power $N_{S_f}:S_f\rightarrow {\bf R} $ of the $L_k^p$-norm and  its $G$-equivariant
 extension are of class  $C^{m_0}$.  Using this   a $C^{m_0}$-smooth cut-off function  supported on $S_f$ can be constructed such 
  that it is still  $C^{m_0}$-smooth viewed in any other local slices.   Hence it gives rise
 a "globally" defined  $C^{m_0}$-smooth cut-off function. By multiplying such a 
 cut-off function $\beta_{S_f}$ with ${\xi_{S_f}}$, we get the desired $C^{m_0}$-smooth
 perturbations for the case of the fixed domain $S^2$.  
 
 \begin{cor}
 	For any  section $\eta_{S_f}$ satisfying $C_1+C_2$, the perturbation 	$\beta_{S_f}{\eta_{S_f}}$   is  of class $C^{m_0}$ viewed in any   local slices. Hence in this sense    it is a    $C^{m_0}$-smooth global section.  
 \end{cor}
 
  Applying this to Gromov/Witten and Floer theories,
  the analytic difficulty in our  global perturbation method \cite{5} is resolved for the case of the fixed domain $S^2$.

 The general  genus zero case is treated in \cite{3, 4}.  However,  the results in this paper alone already essentially contain 
 the "hard" analytic facts needed for  the global perturbation method in \cite{5}. 
 To a large extend,  the proofs in the sequel of this paper \cite{3, 4} are   the   reductions from   the general cases  to  the  case treated 
 in this paper.

This paper is organized as follows. Section 2  gives the proof for Theorem 1.2. Section 3 proves the Theorem 1.1. Section 4
gives a new  proof that 
the evaluation map $E:\Sigma\times  {\widetilde W}(f)\rightarrow M $ define by $E({ x}, g)=g(x)$ is of  class  $C^{m_0}$ as a corollary of Theorem 1.1. This result  was proved in \cite{11} by a different method. Section 5 reproduces the simple  proofs in \cite{13} that   both the $p$-th power $N_{S_f} $ of the $L_k^p$-norm and  its $G$-equivariant extension are of  class $C^{m_0}$.

 In order to concentrate on the main issues related to lack of smoothness of $G$-action, in this paper  the basis facts in the 'usual' analytic foundation of Gromov-Witten/Floer theories are assumed.  
 Other than that, this paper is written only using the  basic facts on Sobolev spaces and standard calculus on Banach  spaces.  They can be found in \cite{7, 2}.

 \section{The  Proof of the Theorem 1.2}

 We make a reduction first.   Fix  an embedding  $i: M\rightarrow {\bf R}^N$. Then we have  the splitting  of    $  E=M\times {\bf R}^N=TM\oplus \nu_M$ of the trivial bundle $E$ into sub bundles,
   where $\nu_M$ is the normal bundle of $M$ of the embedding.
 For a $C^{\infty}$-map $f:\Sigma\rightarrow M$, we get the corresponding   splitting of $L_k^p$-bundles:  $ f^*E\otimes \wedge ^1_\Sigma=(f^*TM\otimes \wedge ^1_\Sigma)\oplus (f^*\nu_M\otimes \wedge ^1_\Sigma)$.     By taking $L_k^p$-sections,  we get the induced splitting  of the Banach bundles ${\cal E}|_{{\tilde W}_{k, p}(f)}\simeq {\cal L}_{{\tilde W}_{k, p}(f)}\oplus {\cal V}_{{\tilde W}_{k, p}(f)}$. Here the fibers at $f$ are the spaces of $L_k^p$-sections of the above corresponding finite dimensional $C^{\infty}$-bundles. There are two parallel transports defined on the trivial bundle
 $E$: the trivial 'flat' one and the one induced from the parallel transports on $TM$ and $\nu_M$. By the usual process in GW-theory, using these two parallel transports,  we get the corresponding local trivializations of ${\cal E}|_{{\tilde W}_{k, p}(f)}$: the 'trivial' one
  ${\cal E}|_{{\tilde W}_{k, p}(f)}\simeq {\tilde W}_{k, p}(f) \times L_k^p(\Sigma, f^*E\otimes \wedge ^1_\Sigma)={\tilde W}_{k, p}(f) \times (\Omega ^1_\Sigma\otimes _{C^{\infty}}L_k^p(\Sigma, {\bf R}^N))$ and the one with respect to the splitting ${\cal E}|_{{\tilde W}_{k, p}(f)}\simeq {\cal L}_{{\tilde W}_{k, p}(f)}\oplus {\cal V}_{{\tilde W}_{k, p}(f)}$.

 \begin{pro}

 	The two trivializations on ${\cal E}|_{{\tilde W}_{k, p}(f)}$ are $C^{\infty}$-equivalent. Moreover, their restrictions to $S_f$ are $C^{\infty}$-equivalent trivializations. 
 	
 \end{pro}
 
 The proof of   this proposition is similar to the   proof of  Proposition 2.2   \cite{13} that identifies the two  local trivializations for ${\cal T}{\cal B}_{k,p}$ of  the tangent bundle  of ${\cal B}_{k,p}$ there.  The details will be given in \cite{17}.
 %Note that restrictions of the trivializations of ${\cal L}_{{\tilde W}_{k, p}(f)}$ and ${\cal V}_{{\tilde W}_{k, p}(f)}$ to $S_f$ are still   the trivializations  obtained by   parallel transport  since .
 
 Now  a section $\xi: {\tilde{\cal B}}_{k, p}\rightarrow {\cal L}_{k, p}$ becomes a section 
 $\xi^E: {\tilde {\cal B}}_{k, p}\rightarrow {\cal E}_{k, p}$.

 The following corollary is an immediate consequence of the proposition above.
 
 \begin{cor}
 	The  inclusion map ${\cal L}_{k, p}\rightarrow {\cal E}_{k, p}$  between Banach bundles is of class $C^{\infty}$. Consequently,	the  section $\xi$ is of class $C^m$ if and only if $\xi^E$ is.
 \end{cor}

 \begin{lemma}
 	Given a section $\gamma:S_f\rightarrow {\cal L}|_{S_f}$, let $\gamma_{{\cal O}_{S_f}}:{\cal O}_{S_f}={\tilde W}|_{k, p}(f)\rightarrow {\cal L}|_{{\cal O}_{S_f}}$ be its $G_e$ euiqvariant extension. Then $(\gamma^E)_{{\cal O}_{S_f}}=(\gamma_{{\cal O}_{S_f}})^E$.
 \end{lemma}

 \proof
 
 For $h\in {\tilde W}(f), $ $(\gamma_{{\cal O}_{S_f}} )^E(h)=i(\gamma_{{\cal O}_{S_f}}(h) )=i(\gamma(h\circ T^{-1}(h))\circ T(h) )=i(\gamma(h\circ T^{-1}(h)))\circ T(h) =\gamma^E(h\circ T^{-1}(h))\circ T(h)=(\gamma^E)_{{\cal O}_{S_f}}(h).$
 
 Here  $i: h^*(TM)\otimes \wedge^1_{\Sigma}\rightarrow  h^*(E)\otimes \wedge^1_{\Sigma}$ is the inclusion map.
 
 Note that $i(s)\circ \phi=i\circ s\circ\phi =i(s\circ\phi)$ for any $\phi \in G_e$ and any section $s: \Sigma \rightarrow h^*(TM)$, which was used  above.

 \QED
 
 Applying this to our case,  we have that  the $G_e$-equivariant extension $(\eta^E_{ S_f})_{{\cal O}_{S_f}}=(\eta_{{\cal O}_{S_f}})^E$.
 
 Thus  $\eta_{{\cal O}_{S_f}}$ is of class $C^{m}$ if and only if $(\eta^E_{ S_f})_{{\cal O}_{S_f}}$ is.   
 
 Hence we only need to prove Theorem 1.2 for a section $\eta=\eta_{S_f}:S_f\rightarrow {\cal E}$.
 
 Now the embedding  $i: M\rightarrow {\bf R}^N$  induces the corresponding $C^{\infty}$-embedding
 ${\bf i}_k:{\tilde {\cal B}}_{k, p} \rightarrow L_k^p=:L_k^p(\Sigma, {\bf R}^N)$ as submanifold. It was proved in Proposition 2.3 of  \cite{13} that the  embedding ${\bf i}_k$ is splitting over an small tubular neighborhood of ${\tilde {\cal B}}_{k, p}$. Consequent  any map $\Psi:{\cal C}\rightarrow {\tilde {\cal B}}_{k, p}$ is smooth if and only if ${\bf i}_{k}\circ \Psi$ is. Now the maps   $\Phi$ and ${\bf i}$ commute each other:  ${\bf i}_{k-m, p}\circ \Phi_{m}^{{\tilde {\cal B}}_{k, p} }=\Phi_{m}\circ {\bf i}_{k, p}.$ 
 
 Let ${\tilde V}(f)$ be the small neighborhood
  of $f$ in $L_k^p$ such that ${\tilde W}(f)={\tilde V}(f)\cap {\tilde {\cal B}}_{k, p}.$  Then  Theorem 1.1  implies the following
  
  \begin{cor}
  	
The corresponding action map 
  	$\Phi^W_m:G_e\times {\tilde W}_{k, p}(f)\rightarrow  {\tilde W}_{k-m, p}(f)$  is of class $C^m$. 
 
 \end{cor}
 
 \begin{theorem}
 	Under the conditions $C_1$ and $C_2$ for $\eta:S_f\rightarrow {\cal E}$, the  section $\eta_{{\cal O}_{S_f}}$ is of class $C^m$. 

 \end{theorem}
 
 \proof

 Under the 'canonical' trivialization, the section $\eta: S_f\rightarrow {\cal E}|_{ S_f}$ becomes $\eta(h)=(h, [\eta](h))$. Here
 $[\eta]:S_f\rightarrow (\Omega ^1_\Sigma)\otimes_{C^{\infty}}L_k^p(\Sigma, {\bf R}^N)$ is the corresponding map.

   Then    $\eta_{{\cal O}_{S_f}}(k)=(\eta(k\circ T^{-1}(k)))\circ T(k)=(k,([\eta](k\circ T^{-1}(k)))\circ T(k) )$. Thus we only need to show that 
 $k\rightarrow ([\eta](k\circ T^{-1}(k)))\circ T(k)$ is of class $C^{m}.$

 Now the map $k\rightarrow ([\eta](k\circ T^{-1}(k)))\circ T(k)$ is decomposed as four maps 
 
 \begin{eqnarray*}
 	k&\rightarrow&(T^{-1}(k), k)\\
 	&\rightarrow&(T(k), k\circ T^{-1}(k))\\
 	&\rightarrow&(T(k), [\eta](k\circ T^{-1}(k)))\\
 	&\rightarrow&([\eta](k\circ T^{-1}(k))\circ T(k).
 \end{eqnarray*}

 The first map $k\rightarrow (T^{-1}(k), k)$ given by $(T^{-1}, Id^W):{\tilde W}_{k, p}(f)={\cal O}_{S_f}\rightarrow G_e\times {\tilde W}_{k, p}(f)$
  is of class $C^{m_0}$; the second $(T^{-1}(k), k)\rightarrow (T(k), k\circ T^{-1}(k))$ given by $IN_{G_e}\times \Phi^W_m :G_e\times{\tilde W}_{k, p}(f)\rightarrow G_e\times({ S_f})_{k-m, p}\subset G_e\times{\tilde W}_{k-m, p}(f)$ is of class $C^{m}$ by the  corollary  above
  , where $IN_{G_e}:G_e\rightarrow G_e$ is the inverse map;  the third 
  $(T(k), k\circ T^{-1}(k))\rightarrow (T(k), [\eta](k\circ T^{-1}(k)))$ given by $Id_{G_e}\times [\eta]_{-m, m}:G_e\times({ S_f})_{k-m, p}(f)\rightarrow G_e\times (\Omega ^1_\Sigma\otimes_{C^{\infty}}L_{k+m}^p(\Sigma, {\bf R}^N))$   is of class $C^{m}$ by the condition $C_1+C_2$ on $\eta$ and the last map $(T(k), [\eta](k\circ T^{-1}(k)))\rightarrow ([\eta](k\circ T^{-1}(k)))\circ T(k)$ given by $\Phi^E_{m}:
 G_e\times (\Omega ^1_\Sigma\otimes_{C^{\infty}}L_{k+m}^p(\Sigma, {\bf R}^N))\rightarrow \Omega ^1_\Sigma\otimes_{C^{\infty}}L_{k}^p(\Sigma, {\bf R}^N)$ is of class $m$ by the Theorem 1.1.

 Now  in the decomposition of $[\eta_{{\cal O}_{S_f}}]$ above, the stated  degree of smoothness for each map is already maximal.
 
 In order to  get the  required smoothness  for $\eta_{{\cal O}_{S_f}}$, the conditions   $C_1$ and $C_2$ on $\eta$  have to be satisfied.
 In this sense, these conditions are necessary.

 \QED

\section {Proof  of   the Theorem  1.1}

%\noindent
%${\bullet }$   ${\bullet }$ Smoothness of the coordinate transformations  of  ${ %{\cal M}}ap_{k, p}$.

First recall the definition of  natural coordinate chart ${\widetilde W}(f)$. For simplicity, assume that $f$ is of class $C^{\infty}.$   It is defined to be $Exp_f:{\hat W}\rightarrow {\widetilde W}(f)=:Exp_f ({\hat W})$, where   ${\hat W}$ consists of all $L_k^p$-sections $\xi$ of the bundle $f^*(TM)\rightarrow \Sigma$ with $\|\xi\|_{k, p}<\epsilon, $ and 
$(Exp_f\xi)(x)=exp_{f(x)}\xi(x)$.
	
 Here the exponential map is taken with respect to an 
$f$-dependent metric defined as follows. Let ${\bf x}=\{x_1,x_2, x_3\}=\{0, 1,\infty\}$ be the three  standard marked points of $S^2$ and ${\bf H}=\{H_1,H_2, H_3\}$ be the three local hypersurfaces of codimension $2$ of $M$ that are transversal to $f$ at the three points $f(x_1), f(x_2)$ and $f(x_3)$. Then  the $f$-dependent metric on $M$ is the one such that each $ H_i, i=1, 2, 3$ is totally geodesic. This can be done by assuming that  the metric near $f(x_i)$ is flat.

Let $h_i\subset T_{f(x_i)}M, i=1, 2, 3$ be the 
tangent spaces of $H_i$ at $f(x_i)$.
Since each $H_i$ is totally geodesic at $f(x_i)$,  $exp_{f(x_i)}(h_i)=H_i$ locally.  Define   ${\hat W}(f, h)$ to  be the  subspace of ${\hat W}(f)$  consisting of  the elements $\eta$  such that $\eta(x_i)\in h_i, i=1, 2, 3.$

Then ${ W}(f, {\bf  H})=:Exp_{f(\bf x)}{\hat W}(f, h)$ is the collection of the  maps $g$ in  ${\widetilde W}(f)$ such that $f({\bf x})\in {\bf H}$. Since the evaluation map $E_{\bf x}:{\widetilde  W}(f )\rightarrow M^3$ for fixed ${\bf x}$ is smooth and obvious transversal to any pint of $M^3$, 
${ W}(f, {\bf  H})=E_{\bf x}^{-1}({\bf H})$ is a $C^{\infty}$ submanifold of ${\widetilde  W}(f )$ with finite co-dimension.
It is a local uniformizer of ${\cal  B}^s$, where ${\cal  B}^s$ is the space of unparametrized stable $L_k^p$-maps. For simplicity, we assume that the elements in  ${\cal  B}^s$ are somewhere injective so that  the action of the reparametrization group on the space ${\widetilde {\cal B}}^s$ of the  parametrized stable $L_k^p$-maps  has local slices given by ${ W}(f, {\bf  H})$(see \cite{16} for a general discussion on stability).

The  following lemma  gives the definition of $T^{-1}:{\tilde W}(f)\rightarrow G_e$  used in Sec. 1. 
It is a consequence of the implicit function theorem applying to 
the evaluation map $E:\Sigma^3\times  {\widetilde W}(f)\rightarrow M^3$ define by $E({\bf x}, g)=(g(x_1), g(x_2), g(x_3))$, assuming the  $C^{m_0}$-smoothness of the evaluation map (proved in
\cite {11} and  next section ).

\begin{lemma}\label{Tmap}
When $\epsilon$ is small enough, there is a $C^{m_0}$-smooth function $T^{-1}:{\widetilde W}(f)={\widetilde W}_{\epsilon}(f)\rightarrow G={\bf PSL}(2, {\bf C})$ such that
for any $k\in {\widetilde W}(f),$ $k\circ T^{-1}(k)\in { W}(f, {\bf  H})$.

\end{lemma}

\proof

Indeed since the evaluation map $E:\Sigma^3\times \{f\}\rightarrow M^3$ is transversal to the local hypersurfaces
 ${\bf H}=\{H_1, H_2, H_3\}$ at $(f(x_1), f(x_2), f(x_3)\in M^3$. By implicit function theorem, the equation $E({\bf y}, k)\in {\bf H}$ has an unique solution ${\bf y}={\bf y}(k)$ for $k\in {\widetilde W}_{\epsilon}(f)$  such that
  ${\bf y}={\bf y}(k)$ is of class $C^{m_0}$ and $\epsilon$-close to ${\bf x}$ when  $\epsilon $ is sufficiently small.
Then there is an unique element $\gamma({\bf y})\in G_e$ smooth in ${\bf y}$ given by ${\bf x}\rightarrow {\bf y}$  so that $(k\circ \gamma ({\bf y}) )({\bf x})\in {\bf H}.$  %Note that the map $k\rightarrow {\bf y}(k)\rightarrow \gamma ({\bf y}(k))$ is of class $C^{m_0}$. 
Define $T^{-1}(k)=\gamma ({\bf y}(k))$.
 \QED

As above, we also denote  the local uniformizer ${ W}(f, {\bf  H})$ by $S_f$. 

 \

  %Here $\phi^*\xi(\gamma, h)=\{\xi (h\circ \gamma)\}\circ (\gamma^{-1})$.

 %Here 
% $$=\Pi(g, f)(\{[\eta]_{f\circ T^{-1}(k)}(k\circ T^{-1}(k))\}\circ T(k))$$
 %$$=:\Pi(g, f)\{(T(k)^*([\eta]_{f}))(k)\}.$$

 %Here the pull-back of $[\eta]_{f}$ by $T(k)$, $T(k)^*([\eta]_{f})$ is defined by 
 %$(T(k)^*([\eta]_{f}))(k)=:\{[\eta]_{f\circ T^{-1}(k)}(k\circ T^{-1}(k))\}\circ T(k).$
 
% Clearly for the smoothness here,  we only need to consider the case  that $g=f$.
  
%To see the degree of the smoothness of $[\eta_{{\cal O}_{S_f}}]_g$, note that 
% $$[\eta_{{\cal O}_{S_f}}]_g(k)$$ $$ =\{\Pi(g, k)\Pi(k , f\circ T(k))\Pi( f\circ T(k), f)\}
 %\{\Pi(f , f\circ T(k))(\{[\eta](k\circ T^{-1}(k))\}\circ T(k))\}.$$

%Now denote the last term $\Pi(f , f\circ T(k))(\{[\eta](k\circ T^{-1}(k))\}\circ T(k))$ by ${\gamma}=:(\Phi\circ (T,[\xi] ))(k),$
%where , and $[\xi]:
%{\widetilde W}(f)\rightarrow {\cal L}_f$ is given by $[\xi](k)=([\eta](k\circ %T^{-1}(k))(T(k))=:T(k)^*([\eta])$.

%The proof of smoothness of $[\eta_{{\cal O}_{S_f}}]_g$  is essentially same as 
 %the one for 

%${\bf \bullet}$ ${\bf \bullet}$  Conditions on $[\eta]$:

%Thus we assume  the following conditions under the assumption that $f$ is $C^{\infty}$ or sufficiently smooth.

%$C_1:$ The map  $[\eta]:S_f\rightarrow  {\cal L}_f$ extends to a $C^{m_0}$-smooth map
%$[\eta]_m:(S_f)_{-m}\rightarrow  {\cal L}_f$ for some $m\geq m_0$

%$C_2:$ The image of  $\eta$ is  lying in $L_{k+m}^p=:({\cal L}_f)_{m}$ with $m\geq m_0$ so that $[\eta]:S_f\rightarrow  ({\cal L}_f)_m$.

We restate Theorem 1.1 in the following  form.

\begin{theorem}\label{Act4}
	The  action map, denoted by  	$\Phi_m:G_e\times L_k^p(\Sigma, {\bf R}^N )\rightarrow  L_{k-m}^p(\Sigma, {\bf R}^N )$ given by $\Phi_m(\phi, h)=h\circ \phi$ is of class $C^m$ assuming that $k-m-2/p>0$.
\end{theorem}

The proof of this theorem will occupy the  rest of this section. 

Clearly we can  replace ${\bf R}^N$ by $ {\bf R}^1$.

We need make another reduction. Let ${\cal V}=\cup_{i\in I}V_i$ and ${\cal U}=\cup_{i\in I}U_i$ be  two  finite open coverings of $\Sigma$ with each $V_i\subset \subset U_i, i\in I $ being   open  disks and $\alpha_i$ and $\beta_i$ being two  partitions  of unit subordinate to ${\cal V}$ and ${\cal U}$ respectively.

Assume that $|I|=l$. Then the maps $I_{\alpha}:L_k^p(\Sigma, {\bf R}^1)\rightarrow \oplus_{i\in I}L_{k, cpt}^p(U_i, {\bf R}^1)\subset \oplus_{i\in I}L_{k}^p(T^2, {\bf R}^1)$ given by $I_{\alpha}(\xi)=(\alpha_1\cdot \xi, \cdots, \alpha_l\cdot \xi)$ and $J_{\beta}:\oplus_{i\in I}L_{k}^p(T^2, {\bf R}^1)\rightarrow L_k^p(\Sigma, {\bf R}^1)$  given by  $J_{\beta}(\eta_1, \cdots, \eta_l)
=\Sigma_{i}^l\beta_i\eta_i$ are bi-linear ( with respect to $(\xi, \alpha)$/ $(\eta, \beta)$) and continuous, and hence are  smooth. Here each $U_i$ is embedded into $T^2$ and ${\beta_i}$ is also considered as a function on $T^2$. Similarly  by considering  each ${\alpha_i}$ is  as a function on $T^2$, $I_{\alpha}$ extends  into  a function  $I_{\alpha}:L_k^p(T^2, {\bf R}^1)\rightarrow \oplus_{i\in I}L_{k, cpt}^p(U_i, {\bf R}^1)\subset \oplus_{i\in I}L_{k}^p(T^2, {\bf R}^1)$. We may assume that
$\beta_i=1$ on  $V'_i$ for some fixed $V_i\subset \subset V'_i\subset \subset U_i$.
Then $J_{\beta}\circ I_{\alpha}\xi=J_{\beta}((\alpha_1\cdot \xi, \cdots, \alpha_l\cdot \xi))=\Sigma_{i}^l\beta_i(\alpha_i\cdot \xi)=\Sigma_{i}^l\alpha_i\cdot \xi=\xi$ is the identity map of $L_k^p(\Sigma, {\bf R}^1)$ as a subspace of  $L_k^p(T^2, {\bf R}^1)$. Here we have used the fact that $\beta_i=1$ on $V_i'.$   

When  $|G_e|$ is small enough, we may assume that  the action $G_e\times V_i\subset V'_i, i\in I.$

 Then for $\xi\in L_k^p(\Sigma, {\bf R}^1))$,  $\Phi (\phi, \xi)=\xi \circ \phi =J_{\beta}\circ I_{\alpha}(\xi \circ \phi )=J_{\beta}(\alpha_1\cdot (\xi \circ \phi ), \cdots, \alpha_l\cdot (\xi \circ \phi ))=J_{\beta}(((\alpha_1 \circ \phi^{-1})\cdot \xi )\circ \phi , \cdots,((\alpha_l \circ \phi^{-1})\cdot \xi )\circ \phi )=J_{\beta}\{ (J_{\alpha\circ \phi^{-1}} (\xi))\circ \phi\}=J_{\beta}\Phi(\phi, J_{\alpha\circ \phi^{-1}} (\xi)).$  Here $\Phi(\phi, J_{\alpha\circ \phi^{-1}} (\xi))$ stands for  $l$-tuple of the actions.
 Since $\alpha$ is of class $C^{\infty}$,   $\alpha =(\alpha_1, \cdots, \alpha_l)\rightarrow \alpha\circ \phi^{-1}$ is 
smooth in $(\alpha, \phi)$ as a map in the relevant spaces. Hence  the map  $(\alpha, \phi, \xi)\rightarrow J_{\alpha\circ \phi^{-1}} (\xi)$ is smooth as well.   Now $J_{\alpha\circ \phi^{-1}} (\xi)$ is a $l$-tuple, each is supported in $U_i$ and hence becomes an element in 
$L_k^p(T^2, {\bf R}^1))$.
The action of $\phi$ in the expression  $\Phi(\phi, J_{\alpha\circ \phi^{-1}} (\xi))$  then is 'supported' inside $U_i$ accordingly for each $i=1, \cdots, l$ so that it can be  extended  into  an action on $T^2$, denoted by ${\hat \phi}_i$ by requiring that ${\hat \phi}_i=id$ outside $V'_i$. Let ${\hat \phi}=({\hat \phi}_1, \cdots,{\hat \phi}_l)$. It is easy to see that the extension ${\hat \phi}$ can be defined in such a way that the map $\phi \rightarrow {\hat \phi}$ is of class $C^{\infty}$ in the corresponding spaces.
Thus  $\Phi(\phi, J_{\alpha\circ \phi^{-1}} (\xi))=\Phi^{T^2}({\hat \phi}, J_{\alpha\circ \phi^{-1}} (\xi)),$
 and we only need to prove  the theorem
  for  $\Phi^{T^2}({\hat \phi}, \eta)$.

Thus we  can indeed  replace the domain $\Sigma$ by  $T^2$ (or ${\bf R}^2$) in the proof of the  Theorem 1.1. 
 The  purpose of this reduction  is to  have  global "coordinate chart"   of the domain.  

In the following proof, we will still  use $\Sigma$ to denote ${\bf R}^2$ or $T^2$.

In the following the  smooth (action) map  is denoted by $T:G_e\times\Sigma\rightarrow \Sigma$   and 
$T_a=T_{\{a\}\times \Sigma}:\Sigma\rightarrow \Sigma$; and   the gradient $(\nabla \eta)=({\partial }_{x_1}\eta , {\partial }_{x_2}\eta)$   of $\eta$ is given using the  global coordinate  $x=(x_1, x_2)$ of $\Sigma.$

Recall that we define  $||\xi||_{k, p}=(\Sigma_{i=0}^k\int_\Sigma |\nabla ^i\xi|^p d\mu_{\Sigma})^{1/p}$ rather than

$\Sigma_{i=0}^k(\int_\Sigma |\nabla ^i\xi|^p d\mu_{\Sigma})^{1/p}.$
\begin{lemma} 
	
For $p>1$, $(\int_0^1 |f(x, t)|dt)^p\leq \int_0^1 |f(x, t)|^pdt$.
Consequently $\|\int_0^1 f (x, t)dt\|^p_{k, p}
\leq  \int_0^1 \|f(x, t)|_{\Sigma\times \{t\}}\|^p_{k, p} dt$.
In particular, $\|\int_0^1\xi\circ T_t dt\|_{k, p}
\leq  C_2(\|T_t\|_{C^k}) \cdot \|\xi\|_{k, p}.$
\end{lemma}

\proof

 The inequality  $(\int_0^1 |f(x, t)|dt)^p\leq \int_0^1 |f(x, t)|^pdt$ follows from  the convexity
of the function $g(x)=x^p$ with $ x>0$ and $p>1.$  

This  implies that $$\|\int_0^1 f (x, t)dt\|_{k, p}^p=\Sigma_{i=0}^k \int_{\Sigma}|\nabla^i_x \int f(x, t)dt |^pd\mu_{\Sigma} \leq \Sigma_{i=0}^k \int_{\Sigma} \int |\nabla^i_x  f(x, t)|^pdt  d\mu_{\Sigma}$$ $$ = \int  \{\Sigma_{i=0}^k \int_{\Sigma} \int |\nabla^i_x  f(x, t)|^p  d\mu_{\Sigma}\}dt
= \int (\|f(x, t)|_{\Sigma\times \{t\}}\|^p_{k, p}) dt .$$

Similarly for a $L_k^p$-function on $\Sigma$ and 
a smooth family of diffeomorphism $T_t:\Sigma\rightarrow \Sigma,$ $$\|\int_0^1\xi\circ T_t dt\|^p_{k, p}$$ $$ = \Sigma_{i=0}^k \int_{\Sigma} |\nabla^i_{\Sigma}(\int_0^1\xi\circ T_t dt)|^p d vol_{\Sigma}$$  can be written, by applying above inequality and the formula of changing variables, as a summation of the terms in the form 

$$ \int_{\Sigma} |\int_0^1(\nabla^i\xi)\circ T_t\cdot F_1({T_t}) dt|^p dvol_{\Sigma}\leq \int_{\Sigma} \int_0^1 |\nabla^i\xi\circ T_t\cdot F_1({T_t})|^p dtdvol_{\Sigma}$$
$$ =\int_0^1\int_{\Sigma}|\nabla^i\xi\cdot F_1({T_t})\cdot  F(Jac_{T_t}) |^p dvol_{\Sigma} dt $$ 
$$\leq  C( \|T_t\|_{C^k})\cdot\|\xi\|_{k, p}^p $$ for some constant
$C( \|T_t\|_{C^k})$ depending on $\|T_t\|_{C^k}$.  Here $F(Jac_{T_t})$ is a  polynomial in $Jac_{T_t}$, $Jac^{-1}_{T_t}$ and their $j$-th derivatives with $j\leq i, $ and similarly for $F_1(T_t)$. 

Then    $\|\int_0^1\xi\circ T_t dt\|^p_{k, p}\leq  C_1 (\|T_t\|_{C^k} )\cdot\|\xi\|_{k,p}^p  $,  hence $$\|\int_0^1\xi\circ T_t dt\|_{L_0}=\|\int_0^1\xi\circ T_t dt\|_{k, p}
\leq  C_2(\|T_t\|_{C^k}) \cdot \|\xi\|_{k, p}=C_3 \|\xi\|_{L_0}.$$  

\QED

Similar proofs to the  last inequality in the lemma above  are  used  throughout  the rest of this section.

In the rest of this section, if there is no confusion, we denote $\Phi_m$ by $\Phi$,  $L_k^p(\Sigma, {\bf R})$ by $L_0$ and $L_{k-m}^p(\Sigma, {\bf R})$ by $L_{-m}$ for short; the norms of the functions  $F(Jac_{T_t}, Jac^{-1}_{T_t})$  in the proofs below will be denoted by $\|Jac_{T_t}||$ or $|| Jac^{-1}_{T_t}\|.$

Let $D_2\Phi$ be the partial derivative of $\Phi$ along $L_0$-direction.
 Since $\Phi(\phi, \xi)=\xi\circ \phi $ is linear in $\xi$,
we have $(D_2\Phi)_{(\phi, \xi)}(\gamma)=\gamma\circ \phi.$
 Hence $D_2:G_e\times {L}_0\rightarrow L( { L}_0,  { L}_{-m})$ is given by $(\phi, \xi)\rightarrow D_2\Phi(\phi, \xi)=\{\Phi_{\phi}=:\Phi(\phi, -):{ L}_0\rightarrow { L}_{-m}\}\in L( { L}_0,  { L}_{-m})$ with
 $\Phi_{\phi}(\eta)=\Phi(\phi, \eta)=\eta\circ \phi.$  Note that $(D_2\Phi){(\phi, \xi)}$ is independent of $\xi$ so that it is factorized as $D_2\Phi=\Psi\circ\pi_1:
 G_e\times { L}_0\rightarrow G_e\rightarrow L(  L_0,  { L}_{-m})$ where 
$\pi_1:G_e\times { L}_0\rightarrow G_e$ is the natural projection map and 
$\Psi:G_e\rightarrow L( {L}_0,  { L}_{-m})$
 is defined by $\Psi(\phi)=\Phi_{\phi}\in L( { L}_0,  { L}_{-m}).$
 Clearly $D_2\Phi$ and $\Psi$ has the same degree of the smoothness.

 Then we have to show  the following proposition.

%${\bullet} $ ${\bullet} $

\begin{pro}
	$\Psi=\Psi_{m-1}: G_e\rightarrow L(L_0,  L_{-m})$ defined by $\Psi(a)=\eta\circ T_a$  for  $\eta\in L_0$ is of class  $C^{m-1}.$  Consequently   $D_2\Phi$ (=$D_2\Phi_{m}$) is  of class  $C^{m-1}$ as well.

\end{pro}

\proof

Let  ${\bf  g}$ be the  'Lie algebra' of $G$. We   identify  $G_e$   with  an open set in ${ \bf g}$. Let $a=(a^1, \cdots, a^l)$ be the coordinate
 of ${\bf g}$ with respect to  a  fixed  basis
 ${\bf a}_j, j=1, \cdots, l.$
 Then the $n$th derivative $D^{n}\Psi:G_e\rightarrow L^n({\bf g}_1\times \cdots\times {\bf  g}_n, L(L_0,  L_{-m})) $ can be identified with  the collection of all $n$-th partial derivatives  at $a\in G_e,$ $\partial ^{\alpha} \Psi(a)$ as  elements in $L(L_0,  L_{-m})$. Here ${\bf g}_i$ is the $i$-th copy of ${\bf g}$,  $L^n({\bf g}_1\times \cdots\times {\bf  g}_n, L(L_0,  L_{-m}))$   is the set of $n$ multi-linear maps from   ${\bf g}_1\times \cdots\times {\bf  g}_n$  to $L(L_0,  L_{-m})$, $\alpha=(\alpha_1, \cdots, \alpha_l)$ with 
 $|\alpha|=|\alpha_1|+\cdots +|\alpha_l|=n$ and $\partial ^{\alpha} \Psi=\partial ^{\alpha_1}_{a_1}\cdots \partial ^{\alpha_l}_{a_l}\Psi. $

 The proposition  is proved in the following steps. 

(I) $\Psi: G_e\rightarrow L(L_0,  L_{-1})$  is continuous. 
 
 \medskip
 \proof   
 
 Indeed  
  $$\|\Psi (a_2)-\Psi(a_1)\|_{L(L_0, L_{-1})} =\max _{\|\eta\|_{L_0}\leq 1}\|\Psi (a_2)(\eta)-\Psi(a_1)(\eta )||_{L_{-1}}$$ 
  
  $$=\max _{\|\eta\|_{L_0}\leq 1}\|\eta\circ T_{a_2}-\eta\circ T_{a_1}||_{L_{-1}} =\max _{\|\eta\|_{L_0}\leq 1}\|\int_0^1 {\frac{d}{dt  } }(\eta\circ T_{a_1+t(a_2-a_1)})dt ||_{L_{-1}} $$
  
  $$=\max _{\|\eta\|_{L_0}\leq 1}\|\int_0^1 (\nabla\eta\circ T_{a_1+t(a_2-a_1)})\cdot (D_tT_{a_1+t(a_2-a_1)}(a_2-a_1)) dt ||_{L_{-1}} $$  $$ \leq \max _{\|\eta\|_{L_0}\leq 1}\int_0^1\| (\nabla\eta\circ T_{a_1+t(a_2-a_1)})\cdot (D_tT_{a_1+t(a_2-a_1)}(a_2-a_1))||_{L_{-1}} dt$$
  $$\leq \max _{\|\eta\|_{L_0}\leq 1}\int_0^1\| \nabla\eta\|_{L_{-1}} \cdot || Jac_{T_{a_1+t(a_2-a_1)}}\|_{C^{k-1}}\cdot \|D_tT\|_{C^{k-1}}\cdot||a_2-a_1|| dt$$
  $$\leq \int_0^1 || Jac_{T_{a_1+t(a_2-a_1)}}\|_{C^{k-1}}\cdot \|D_tT\|_{C^{k-1}} dt\cdot||a_2-a_1|| =C||a_2-a_1||$$ 
  
  for some constant only depending on $T$. 
  \QED
  
 In fact we have proved that $\Psi_1$ above is Lipschitz, though we do not need this.

 (II)  The map  $\Phi^D: G_e\rightarrow L(L_0,  L_{-1})\oplus  L(L_0,  L_{-1})$  given by $\Phi^D(a)(\eta) =\nabla\eta\circ (T_a)$ exists for $\eta \in L_0$. % Here $T:G_e\times\Sigma\rightarrow \Sigma$  is the smooth (action) map and 
  %$T_a=T_{\{a\}\times \Sigma}:\Sigma\rightarrow \Sigma$; $(\nabla \eta)=({\partial }_{x_1}\eta , {\partial }_{x_2}\eta)$ is the gradient of $\eta$, where $x=(x_1, x_2)$ is the  global coordinate of $\Sigma.$ 
  Moreover, $\Phi^D: G_e\rightarrow L(L_0,  L_{-2})\oplus  L(L_0,  L_{-2})$
 is continuous.

\medskip
\noindent
Note that $\eta\rightarrow \nabla \eta= (\partial_{x_1}\eta, \partial_{x_2}\eta )$ is considered as a map: $L_0\rightarrow L_{-1}\oplus L_{-1}.$

\medskip
\proof

 Existence part  follows from the definition.  Indeed denote $L(L_0,  L_{-1})\oplus  L(L_0,  L_{-1})$ by $\oplus ^2L(L_0,  L_{-1})$ for short.
 
 Then for any $a\in G_e,$ $$\|\Phi^D(a)\|_{\oplus ^2L(L_0,  L_{-1})}=Max_{||\eta||_{L_0}\leq 1}\|\Phi^D(a)(\eta)\|_{ L_{-1}\oplus  L_{-1}}$$ 
 $$=Max_{||\eta||_{L_0}\leq 1}\|\nabla\eta\circ (T_a)\|_{ L_{-1}\oplus  L_{-1}}\leq  ||Jac^{-1}_{T_a}\|_{C^{k}}\cdot Max_{||\eta||_{L_0}\leq 1}\|\nabla\eta\|_{L_{-1}\oplus  L_{-1}}$$
 $$\leq ||Jac^{-1}_{T_a}\|_{C^{k}}.$$

 The continuity for each summand then follows from the continuity  in (I) by shifting the index by one and replacing $\eta$ by $\nabla \eta$.
 
\QED 

 (III) The partial derivatives  $\partial _{j} \Psi=:\partial _{a_j} \Psi:  G_e\rightarrow     L(L_0,  L_{-2})            $   exists  and $\partial _{j} \Psi (a)(\eta)=   \{ \nabla\eta\circ (T_a) \}\cdot (\partial _jT)_a =\{\Phi^D(a)(\eta)\}\cdot (\partial _jT)_a  $.  In other words, upto a harmless factor, $\partial _{j} \Psi $ is essentially  just
 $\Phi^D$ so that it is continuous by  (II) above.
 
 \medskip
 \proof
 
 Indeed $$||\Psi(a+t{\bf a}_j)-\Psi(a)-\partial _{j} \Psi (a)(t{\bf a}_j))
\|_{L(L_0, L_{-2})} $$
 
  $$=\max _{\|\eta\|_{L_0}\leq 1}\|\eta \circ T_{a+t{\bf a}_j}-\eta \circ T_{a}-t\cdot \{ \nabla\eta\circ (T_a) \}\cdot (\partial _jT)_a||_{L_{-2}}$$
 
 $$=\max _{\|\eta\|_{L_0}\leq 1}\|\int_0^1\{{\frac{d}{ds  } }(\eta \circ T_{a+st{\bf a}_j})-t\cdot \{ \nabla\eta\circ (T_a) \}\cdot (\partial _jT)_a\} ds||_{L_{-2}}$$
 
 $$=\max _{\|\eta\|_{L_0}\leq 1}\|\int_0^1 t\{\nabla\eta \circ T_{a+st{\bf a}_j}\}\cdot (\partial _jT)_{a+st{\bf a}_j}-t\cdot \{ \nabla\eta\circ (T_a) \}\cdot (\partial _jT)_a\} ds||_{L_{-2}}$$
 
 $$\leq |t| \max _{\|\eta\|_{L_0}\leq 1}\int_0^1\|\{\nabla\eta \circ T_{a+st{\bf a}_j}\}\cdot (\partial _jT)_{a+st{\bf a}_j}- \{ \nabla\eta\circ (T_a) \}\cdot (\partial _jT)_a ||_{L_{-2}}ds$$

 $$= |t| \max _{\|\eta\|_{L_0}\leq 1}\int_0^1\|\int_0^1{\frac{d}{d\nu   } }\{\nabla\eta \circ T_{a+\nu st{\bf a}_j}\cdot (\partial _jT)_{a+\nu st{\bf a}_j}\}d\nu||_{L_{-2}}ds$$
 
  $$\leq |t| \max _{\|\eta\|_{L_0}\leq 1}\int_0^1\int_0^1\|{\frac{d}{d\nu   } }\{\nabla\eta \circ T_{a+\nu st{\bf a}_j}\cdot (\partial _jT)_{a+\nu st{\bf a}_j}\}||_{L_{-2}}d\nu ds$$

 $$\leq |t| ^2\max _{\|\eta\|_{L_0}\leq 1}\int_0^1\int_0^1\|\nabla^2\eta \circ T_{a+\nu st{\bf a}_j}||_{L_{-2}}\cdot ||(\partial _jT)_{a+\nu st{\bf a}_j}||^2_{{C^{k-2}}}$$ $$+\|\nabla\eta \circ T_{a+\nu st{\bf a}_j}||_{L_{-2}}\cdot ||(\partial _jT)_{a+\nu st{\bf a}_j}||_{{C^{k-1}}}  d\nu ds$$
 
 $$\leq |t| ^2\max _{\|\eta\|_{L_0}\leq 1}\int_0^1\int_0^1\|\eta \|_{L_{0}}\cdot  \|Jac^{-1}_{T_{a+\nu st{\bf a}_j}}||_{C^{k-2}}\cdot (1+||T_{a+\nu st{\bf a}_j}||_{C^{k}})^2d\nu ds$$
 
$$=|t| ^2\cdot  \int_0^1\int_0^1 \|Jac^{-1}_{T_{a+\nu st{\bf a}_j}}||_{C^{k-2}}\cdot (1+||T_{a+\nu st{\bf a}_j}||_{{C^{k}}})^2d\nu ds.$$

%Here  as well as in (I), we have used $||\eta \circ T_a||_{L_0}\leq \|\eta\|_{L_0}\cdot \|Jac^{-1}_{T_a}\|_{C^k}$ in  the lemma above. 

  \QED

  (IV)  Inductively, we assume that the $(m-1)$-th partial derivatives with $|\alpha|=m-1$ and $\alpha=(\alpha_1, \cdots, \alpha_l)$ is a  map $\partial^{|\alpha|}_{\alpha}\Psi=:\partial^{\alpha_1}_{j_1}\cdots \partial^{\alpha_l}_{j_l}\Psi:G_e\rightarrow  L(L_0,L_{-m})$  and that for $a\in G_e$ and $\eta\in L_0$,  $\partial^{|\alpha|}_{\alpha}\Psi(a) (\eta)$ is a summation of the terms $\partial^{\gamma}\eta\circ T_a=:\partial_{x_1}^{\gamma_1}\partial_{x_2}^{\gamma_2}\eta$ multiplying with $(\partial^{\beta}T)_a=(\partial^{\beta_1}_{j_1}\cdots \partial^{\beta_l}_{j_l}T)_a$
  with $\partial^{\beta}T$ along $G_e$-directions. Here $|\gamma|$ and $|\beta|$ is less than or equal to  $m-1$. We have to show that 
  
  (A) $\partial^{|\alpha|}_{\alpha}\Psi:G_e\rightarrow  L(L_0,L_{-m})$ with $|\alpha|=m-1$   is continuous.
  
  (B) $\partial_j\partial^{|\alpha|}_{\alpha}\Psi:G_e\rightarrow  L(L_0,L_{-m-1})$ with $|\alpha|=m-1$ is a summation of the terms $\partial^{\gamma}\eta\circ T_a$ multiplying with $(\partial^{\beta}T)_a$
  with $\partial^{\beta}T$ along $G_e$-directions. Here $|\gamma|$ and $|\beta|$ is less than or equal to  $m$.  The  proof of A and B will complete the induction.
  
  To prove A and B, note that the factor $(\partial^{\beta}T)_a$ do not affect anything here.  Denote the term $\partial^{\gamma}\eta\circ T_a$  by  $\Phi^{D, \gamma}(a)
(\eta)$ . Then  $  \Phi^{D, \gamma}:G_e\rightarrow L(L_0, L_{-|\gamma|})$. Note that $\Phi^{D, \gamma}$  has essentially the same type of form  as $\Phi^D$ has in (II) above.
By (II) with an obvious degree shifting, (i) the image of $\Phi^{D, \gamma}$ is indeed in $L(L_0, L_{-|\gamma|})$ and (ii ) $\Phi^{D, \gamma}:G_e\rightarrow L(L_0, L_{-|\gamma|-1})$ is continuous.  
This proves A by letting $|\gamma|=m-1$.

 Similarly,  $  \Phi^{D, \gamma}$ and $\Psi$ are also the "same type". Hence if we shift the  degree of the target above by one  further and consider $\partial_j  \Phi^{D, \gamma}:G_e\rightarrow L(L_0, L_{-m-1})$.  The  essentially the same argument of  III shows that B is true for $\partial_j  \Phi^{D, \gamma}$  for $|\gamma|=m-1$.  This implies that $B$ is true for $\partial_j  \partial^{|\alpha|}_{\alpha}\Psi$ by the definition of $\Phi^{D, \gamma}$ above.

  \QED

 Next  consider the partial derivative $D_1\Phi=D_G\Phi$ along $G_{e}$-direction
 for $\Phi_n:G_e\times {L}_0\rightarrow { L}_{-n}$. 
 Since $G_e$ is finite dimensional so that   $L({\bf g}, L_{-n})$ is a finite sum of $L_{-n}$,  by the discussion before, the partial derivative $D_1\Phi_n:G_e\times L_0\rightarrow L({\bf g }, L_{-n})$ can be identified with the collection of ${\partial}_{a_j}\Phi_n=:{\partial}_{j}\Phi_n: G_e\times L_0\rightarrow L_{-n} , j=1, \cdots, l.$ Denote ${\partial}_{a_j}\Phi_k: G_e\times L_0\rightarrow L_{-k} , j=1, \cdots, l, $ by $({\partial}_{j}\Phi)_k$.

 We  have to show the following lemma.
 
% ${\bf \bullet}$ ${\bf \bullet}$ 
 
 \begin{lemma}
  For each $j$, $({\partial}_{j}\Phi)_n:G_e\times L_{0}\rightarrow L_{-n}$ is of class $C^{n-1}.$
\end{lemma}

 \proof
 
 We claim that  as a map $({\partial}_{j}\Phi)=:({\partial}_{j}\Phi)_0:G_e\times L_0\rightarrow L_{-1}$,  ${\partial}_{j}\Phi(a, \eta)=(\nabla \eta)\circ (T_{a})
 \cdot {\partial }_{{\bf a}_j}(T)_{a}.$
 
 Indeed for $(a, \eta)\in G_e\times L_0$,  $$\|\Phi(a+t{\bf a}_j, \eta)-\Phi(a, \eta)-t\cdot  (\nabla \eta)\circ (T_{a})
 \cdot {\partial }_{{\bf a}_j}(T)_{a}\|_{k-1, p }$$ $$ =|t|\cdot \| \int_0^1 \{(\nabla \eta)\circ (T_{a+\nu t {\bf a}_j})\}
 \cdot {\partial }_{{\bf a}_j}(T)_{a+\nu t{\bf a}_j}\}-\{(\nabla \eta)\circ (T_{a})
 \cdot {\partial }_{{\bf a}_j}(T)_{a}\}d\nu \|_{k-1, p }$$ $$\leq |t|\cdot \int_0^1 \| \{(\nabla \eta)\circ (T_{a+\nu t {\bf a}_j})\}
 \cdot {\partial }_{{\bf a}_j}(T)_{a+\nu t{\bf a}_j}\}-\{(\nabla \eta)\circ (T_{a})
 \cdot {\partial }_{{\bf a}_j}(T)_{a}\}d\nu \|_{k-1, p }.$$
 
 Now for fixed $\eta$ with $(\nabla \eta)\in L_{-1}=L_{k-1}^p$, the continuity of the map $\Phi: G_e\times L_{-1}\rightarrow L_{-1} $ given $(\xi, a)\rightarrow  \xi\circ T_a$ implies that the integrand above goes to zero as 
 $t\rightarrow 0$ uniformly for $\nu$ since $|\nu \cdot t|\leq |t|.$ This implies that the integral above goes to zero as  $t\rightarrow 0$,  and hence the  claim is proved.

 Thus up a harmless factor ${\partial }_{{\bf a}_j}(T)_{a}$ ,   ${\partial}_{a_j}  \Phi $  is essentially given by   $\Phi_j:G_e\times L_0\rightarrow L_{-1}$ defined by $\Phi_j (a, \eta)=(\nabla \eta)\circ (T_{a}).$ Now $\Phi_j$ is of the same type as the original $\Phi:G_e\times L_0\rightarrow L_{0}$ with  a degree shifting by  $1$ and replacement  of $\eta$ by $\nabla \eta$. Hence we are in the position to apply induction to conclude that
 $\Phi_j:G_e\times L_{0}\rightarrow L_{-n}$ is of class $C^{n-1}$ so that  $({\partial}_{j}\Phi)_n:G_e\times L_{0}\rightarrow L_{-n}$ is of class $C^{n-1}.$

 This finishes the proof of the main theorem of this section under the assumption  of the continuity    of  $ (\Phi)_0$.
 
 \QED

 \begin{pro}\label{ContAct}
 $(\Phi)_0:G_e\times L_0\rightarrow L_{0}$ is continuous.
 \end{pro}
 
 \proof

 $$ \|\Phi(a+t, \eta+\gamma)- \Phi(a, \eta)\|_{L_0}$$ $$\leq 
 \|\Phi(a+t, \eta+\gamma)- \Phi(a+t, \eta)\|_{L_0}+\|\Phi(a+t, \eta)- \Phi(a, \eta)\|_{L_0}=A+B.$$ 
 
 Now $$A=\|(\eta+\gamma)\circ T_{a+t}-\eta\circ T_{a+t}\|_{L_0}$$ $$ =\|\gamma\circ T_{a+t}\|_{k, p}\leq C(J(T), \|T\|_{C^k})\cdot \|\gamma\|_{k, p}$$ for a constant  $C(J(T), \|T\|_{C^k})$ depending only on $J(T)$ and $ \|T\|_{C^k}$. Here and below  $J(T)$  is the Jacobian of $T$.
 Hence $A<\epsilon/8$ if we choose $|\gamma\|_{L_0}< \delta_1=\epsilon/(8\cdot C(J(T) \|T\|_{C^k}).$
 
 To estimate $B$,  choose  a smooth $\xi$ such that $\|\xi-\eta \|_{k, p}< \delta_2.$
  Then $$B\leq \|\Phi(a+t, \xi)- \Phi(a, \xi)\|_{L_0}+\|\Phi(a+t, \eta)- \Phi(a+t, \xi)\|_{L_0}+\|\Phi(a, \eta)- \Phi(a, \xi)\|_{L_0}$$ $$ =B_1+B_2+B_3.$$
  
  Clear the estimate for $A$ above is applicable to both $B_2$ and $B_3$  so that for $i=2, 3,$
  $B_i \leq C(J(T), \|T\|_{C^k})\cdot \|\xi-\eta\|_{k, p}\leq C(J(T), \|T\|_{C^k})\cdot \delta_2 < \epsilon/8$ if we choose $\delta_2=\epsilon/(8\cdot C(J(T) \|T\|_{C^k})).$
  
  Now $B_1=\|\xi\circ T_{a+t}-\xi\circ T_{a}\|_{k, p}=\|\int_0^1{\frac {d}{d\nu}}\{\xi\circ T_{a+\nu t}\}d\nu\|_{k, p}$. 
  	 Applying the estimate in $I$  or $III$ above (or its argument) to the case here, we get $B_1\leq C(J(T), \|T\|_{C^{k+1}})\cdot \|t\|\cdot \|\xi||_{k+1, p}<\epsilon/8$ if we choose $|t|<\delta_3=\epsilon/(8\cdot C(J(T) ,\|T\|_{C^{k+1}})\cdot  \|\xi||_{k+1, p}).$

 \QED
 
% We now record a  more precise estimate for term $A$ in above proof  which will be used in the sequel of this paper.
 
% \begin{lemma}
% For  any $L_k^p$-maps $h_1$ and $h_2$, $\|h_1\circ T_t- h_2\circ T_t ||^p_{k, p}<C(\|T_t\|_{C^k})\|h_1-h_2\|^p_{k, p}$, where $C(\|T_t\|_{C^k})$ is a fixed continuous  function of $\|T_t\|_{C^k}$. In particular, $lim_{t\rightarrow t_0}C(\|T_{t_0}\|_{C^k})=C(\|T_t\|_{C^k}).$
% Consequently, the action map $\Phi:G_e\times  L_k^p(\Sigma, M)\rightarrow  L_k^p(\Sigma, M)$ given by 
% $\Phi(t,h)=h\circ T_t$ is uniformly continuous for fixed but bounded $t$ in spite of the fact that $\Phi$ is not uniformly continuous. More specifically, For  fixed  $C$ with $|t|\leq C, $ given $\epsilon>0, $ there is a $\delta$ such that if $\|h_1-h_2\|_{k, p}<\delta, $  $\|\Phi(t, h_1)-\Phi(t, h_2)|_{k, p}<\epsilon$.  
 
 %\end{lemma}
 
% The lemma is proved by using the following estimate  repeatedly.
 
% $$\int_{\Sigma}||h\circ T_t (x)||^p dx=\int_{\Sigma}||h(y)||^p\cdot Jac^{-1}(T_t)(y) dy.$$
 
 \noindent
 
 \section { Smoothness of the evaluation  map}

 The following theorem was proved in \cite {11} by a different method. Here we reprove it as a corollary of  the above Theorem 1.1.

 \begin{theorem}
 	Assume that $p>2$.	The $m$-fold total evaluation map:$E_m:{\widetilde {\cal B} }\times \Sigma^m\rightarrow M^m$  is of class $C^{m_0}$,  where ${\widetilde {\cal B} }={\widetilde {\cal B} }_{k, p}$ is the collection of $L_k^p$ maps from $\Sigma\simeq S^2$ to $M$.
 \end{theorem}
 
 \proof
 
 We make some reductions.  (i ) It is sufficient to prove the theorem for $E_1$ and  $ {\widetilde W }(f)$; (ii) It is reduced further to show  that for any $x_0\in \Sigma$, the map $E_{D(x_0)}=E_1|_{{\widetilde W }(f)\times D(x_0)}:{\widetilde W }(f)\times D(x_0)\rightarrow M$ is of class $C^{m_0}$, where $D(x_0)$ is the small disc on $\Sigma$ centered at $x_0.$

  Let $\phi=:\phi_{D(x_0)}:D(x_0)\times\Sigma\rightarrow \Sigma$ be the smooth "action" map
 such that for any $x\in D(x_0)$, the restriction $\phi_x: \{x\}\times \Sigma\rightarrow \Sigma $ is a diffeomorphism with the property that (a) $\phi_x:D(x_0)\rightarrow   D(x)$ by  a "translation", in particular $\phi_x(x_0)=x$ and (b) $\phi_x$ is  the identity map outside a larger disc centered at $x_0$.

 Denote $D(x_0)$ by $G_e$. Consider the corresponding action map induced by $\phi$ above, $\Phi_{k-1}:G_e\times {\widetilde W }(f)=D(x_0)\times {\widetilde W }(f)\rightarrow ({\widetilde W }(f))_{-(k-1)}$. Then by  Theorem 1.1 and Corollary 2.1 $\Phi_{k-1}$ is of class
 $C^{m_0}$ since $m_0=[k-2/p]=k-1.$

 Now $ ({\widetilde W }(f))_{-(k-1)}$ consists of $L_1^p$-maps that are continuous by  the assumption  $p>2$,  so that the evaluation map at a fixed point $x_0, $ $E_{x_0}:({\widetilde W }(f))_{-(k-1)}\rightarrow M$ is well-defined. Clearly $E_{x_0}$ is smooth of class $C^{\infty}$:  in the local charts $Exp_f:
 {\hat W}(f)\rightarrow {\widetilde W }(f)$	and $exp_{f(x_0)}:B(f(x_0))\rightarrow M$,   the map $E_{x_0}$ is corresponding to
 ${\hat E}_{x_0}$ given by ${\hat E}_{x_0}(\xi)=\xi(x_0)$, which is linear.
 Here  $B(f(x_0))\subset T_{f(x_0)}M$  is the corresponding  ball.
 
 Recall that $\Phi_{k-1}(x, g)=g\circ \phi_x.$
 Now  $E_{D(x_0)}=E_{x_0}\circ \Phi_{k-1}:G_e\times  {\widetilde W }(f)\rightarrow M$, hence is of class $C^{m_0}$. Indeed for any $(x, g)\in D(x_0)\times  {\widetilde W }(f), $
 $E_{x_0}\circ \Phi_{k-1}(x, g)=E_{x_0}(g\circ \phi_x)=g(\phi_x(x_0))=g(x)=E_1( g, x)=E_{D(x_0)}( g, x).$
 
 \QED

 \section{Smoothness of cut-off functions} 
 
 We first prove  the smoothness of the $p$-th power $N_{k}$ of $L_k^p$-norm on ${\tilde B}={\tilde B}_{k,p}$.
It is sufficient to consider the case  ${\tilde B}=L_k^p(\Sigma, {\bf R}^M)$.

\begin{theorem}
	
 Assume that $p$ is an positive even integer. Then $N_{k}:L_k^p(\Sigma, {\bf R}^M)\rightarrow {\bf R}$ defined  by $N_k(f)=\|f\|_{k, p}^p=
\Sigma_{i=0}^k \|\nabla^if\|_{p}^p$ is smooth.

\end{theorem}

\proof

Clearly  it is sufficient to prove the case that $k=0$.

Then $N_0(f)=\int_{\Sigma}<f,  f>^m \, dvol_{\Sigma}$, where $m=p/2$.
Since $I:L^1(\Sigma, {\bf R})\rightarrow {\bf R}^1$  given by $I(h)=\int_{\Sigma}h \, dvol_{\Sigma}$  is linear and hence smooth, we only need to show that
the map $P:L^p(\Sigma, {\bf R}^M)\rightarrow L^1(\Sigma, {\bf R})$ given by 
$P(f)=<f,  f>^m $ is well-defined and smooth.
 Indeed, $P(f)=M\circ \Delta (f)$ where  $\Delta=\Delta_p:L^p(\Sigma, {\bf R}^M)\rightarrow (L^p(\Sigma, {\bf R}^M))^p$ is the $p$-fold diagonal map given by 
 $\Delta(f)=(f_1, \cdots, f_p)$ with $f_1=\cdots =f_p=f,$ and $M:(L^p(\Sigma, {\bf R}^M))^p\rightarrow L^1(\Sigma, {\bf R})$ given by $M(f_1, \cdots, f_{p})
=<f_1, f_2>
\cdots <f_{p-1}, f_p>.$   Clearly $\Delta$ is smooth. We have to show that $M$ is well-defined and continuous. Indeed, by Holder  inequality,
$$\|M(f_1, \cdots f_p)\|_{L^1}=||<f_1, f_2>
\cdots <f_{p-1}, f_p>||_{L^1}\leq ||\, |f_1|\cdot |f_2|
\cdots |f_p|\, ||_{L^1}$$ $$ \leq \|f_1\|_p\cdots \|f_p\|_p.$$ Hence $M$ is well-defined, continuous and multi-linear. This implies that $M$ is smooth.

\QED

\medskip
\noindent
\textbf{Note:}
The theorem above  was proved in \cite{8, 9, 12} by a different method.  This  result was mentioned without proof in various articles or books. It was used in Sec. 4 of \cite{5}.  Until recently  the author   found that its first formal proof  was  given  in \cite{10} as a consequence of the general results there.  Comparing with all those proofs, the proof above is much simpler.

Now consider the action map $\phi:G_e\times \Sigma\rightarrow \Sigma\subset  {\bf R}^{S}$.  Here we have fixed an embedding of $\Sigma$   into ${\bf R}^{S}.$  We require that $\phi$ is of class $C^{\infty}$  as a  map  into ${\bf R}^S$.% In particular for any $k$, $\|\psi\|_{C^k}$ is bounded. So is 
%$\|Jac(\psi_g)\|_{C^k}$. Here the $C^k$-norm is taken over $G_e\times \Sigma$. 

\begin{lemma}
For any $k$, 	the following maps are of class $C^{\infty}$: $Jac: G_e\rightarrow C^k(\Sigma, {\bf R})$  defined  by $Jac (g)=Jac(\phi_g)$,   $\phi_G:G_e \rightarrow C^k(\Sigma, {\bf R}^{S})$ defined by $\phi_G(g)=\phi_g$ and $D^l\phi_G:G_e \rightarrow C^k(\Sigma, {\bf R}^{S})$ defined by $(D^l\phi_G)(g)=D^l\phi_g$. Here $D^l\phi_g$ is the  partial derivatives of  $\phi_g$ along $\Sigma$.

\end{lemma}

\proof 

We leave the direct  proof  to the readers. It is straightforward but somewhat tedious. 

For an indirect proof, note that $\phi_G(g)=\phi_g=\Phi (id_{\Sigma}, g)$. Since $id_{\Sigma}:\Sigma\rightarrow {\bf R}^S$ is of class $C^{\infty}$, the conclusion follows from the discussion before on the action map $\Phi $.

\QED

Denote the  induced action of $\phi$  by $\Phi:G_e\times L_k^p(\Sigma, {\bf R}^M)\rightarrow L_k^p(\Sigma, {\bf R}^M)$ given by $\Phi(g, f)=f\circ \phi_g.$ Then the function $F_k:G_e\times L_k^p(\Sigma, {\bf R}^M)\rightarrow {\bf R}$ is defined by $F_k=N_k\circ \Phi$.

\begin{theorem}
Assume that $p$ is an positive even integer, and that all $\phi_g$ are orientation-preserving  diffeomorphisms.
Then the function $F_k$ is  smooth.

\end{theorem}

\proof

 Since $F_k(g, f)=\Sigma_{i=0}^k\int_{\Sigma}|D^i(f\circ \phi_g)|^pdvol_{\Sigma},$
 it  is  a summation of the terms in the form $$\int_{\Sigma}|(D^jf)\circ \phi_g|^p\cdot | D^l\phi_g|^p
 dvol_{\Sigma}$$ $$ =\int_{\Sigma}(|(D^jf)|\cdot | D^l\phi_G(g^{-1})|)^p\cdot  (Jac^{-1} (g))dvol_{\Sigma}$$ $$ =I((Jac^{-1} (g))\cdot P((D^jf)\cdot D^l\phi_G(g^{-1}))$$ with $j,l \leq k.$

  Thus each such term is  decomposed as the following smooth maps:
 (i) First map is $(f,g )\rightarrow |(D^jf)|^2\cdot | D^l\phi_G(g^{-1})|^2$ from $L_k^p\times G_e\rightarrow L^p $.  (ii) The next is $|(D^jf)|\cdot |D^l\phi_G(g)|\rightarrow P((D^jf)\cdot D^l\psi_G(g))$ given by  the   map $P:L_k^p\rightarrow L_1$. (iii)
 The third map $P((D^jf)\cdot D^l\psi_G(g))\rightarrow (Jac^{-1} (g))\cdot P((D^jf)\cdot D^l\psi_G(g)$  is from  $C^{k}\times L^1\rightarrow  L^1$ given by $(h, f)\rightarrow h\cdot f$.  (iv) The  last map is  the integration  $I:L^1\rightarrow {\bf R}$. 
 
 The conclusion follows.

\QED

An immediate corollary is the following.

\begin{cor}
Let $(N_{k})_{S_f}$  be the restrictions of $N_{k}$  to the  $C^{m_0}$ local slices $S_f$, then its $G_e$-equivariant
extension is of class $C^{m_0}$.
\end{cor}


\begin{thebibliography}{0}
 	
 	%\bibitem{1}Atiyah, M. F.  and Hirzebruch, F:  \textsl{Vector bundles and homogeneous spaces},
 	%	Proceedings of Symposiam of the American Mathematical Society, Vol. III (1960), pp.7-38.
 	\bibitem{8} B. Chen, and A. Li, and B. Wang:  \textsl{Virtual Neighborhood Technique for Pseudo-holomorphic Spheres}, 	arXiv: math.1306.3276[math.GT], 2013.
 	
 	\bibitem{9} X. Chen:  \textsl{Note on Partition of Unit}, 	 Electronic file, Feb. 8, 2013.
 	
 	\bibitem{10} R. Bonic and J.Frampton:  \textsl{Smooth functions on  Banach Manifolds}, J. Math. Mech.J.  Vol. 2, 1966.
 
 	\bibitem{1} Hofer, H:  \textsl{A general Fredholm Theory and applications}, 	arXiv: math.0509366[math.SG], 2005.
 	
 	
 	\bibitem{2}Lang, S:  \textsl{Differential manifolds}, Springer-Verlag, 1972.
 	
 	
 	\bibitem{3} Liu, G :  \textsl{ Higher-degree smoothness of perturbations II}, Preprint, arXiv:1809.08620  [math.SG], 2018.
 	
 	
 	\bibitem{4}  Liu, G :  \textsl{Higher-degree smoothness of perturbations III}, Preprint arXiv:1810.00459 [math.SG], 2018.
 	
 	\bibitem{11}  Liu, G :  \textsl{$C^{m_0}$-smoothness of the evaluation maps}, Preprint,  arXiv:1312.216  [math.SG], 2013.
 	
 	\bibitem{12}  Liu, G :  \textsl{Weak  Smoothness in GW Theory}, Preprint,  arXiv:1310.7209  [math.SG], 2013.
 
 \bibitem{13}  Liu, G :  \textsl{Geometric equivariant extension of sections in GW theory}, Preprint,  arXiv:arXiv:1509.06722  [math.SG], 2015.
 
 \bibitem{14}  Liu, G :  \textsl{ $C^l$-Smoothness of the Action Map
 	$\Psi:G \times  L^p_k(\Sigma, M)\rightarrow  L^p_{k-l}(\Sigma, M)$}, Preprint, 2014.
 
 %\bibitem{15}  Liu, G: \textsl{ $C^1$-Smoothness of  Perturbations},  In Preparation. %Preprint,  arXiv:arXiv:1509.06722  [math.SG], 2019.
 
 \bibitem{16}  Liu, G: \textsl{ On stability of nodal $L_k^p$-maps I }, Preprint,  arXiv:1509.07187  [math.SG],
2015.
 
 \bibitem{17}  Liu, G: \textsl{A functorial framework for analysis on mapping spaces  }, In Preparation.	
 	

 	\bibitem{5} Liu, G. and  Tian, G:  \textsl{ Floer homology and Arnold conjecture }, J. Differential Geom. 49 (1998), no. 2. 
 	
 	
 	\bibitem{6} McDuff, D and   Wehrheim, D. McDuff and   Wehrheim,K :  \textsl{Smooth Kuranishi     Structures with trivial isotropy},Preprint 2012. 
 	
 	\bibitem{7} Palais, R:   \textsl{ Foundations of Global Non-linear Analysis},  Benjamin, INC, 1968.
 	
 	
 	
 	
 	
 	
 	
 	
 	
 	
 	
 \end{thebibliography}
\end{document}